\documentclass[10pt,a4paper]{article}
\usepackage{amsmath}
\usepackage{amsfonts}
\usepackage{amssymb}
\usepackage[normalem]{ulem}
\usepackage{enumitem}
\usepackage[round]{natbib}

\usepackage{float}
\floatstyle{boxed} 
\restylefloat{figure}

\title{Kurt G\"odel's reception of Charles Hartshorne's ontological proof\footnote{This is an AAM before peer-review. The final version is edited and restructured from this version and published in:  E. Ramharter (ed.), The Vienna Circle and Religion. Vienna Circle Institute Yearbook. Springer (2021). }}
\author{ Annika Kanckos\footnote{
		annika.kanckos@helsinki.fi}, and Tim Lethen\footnote{lethen@cs.uni-saarland.de} }

\begin{document}
	\maketitle
	\begin{abstract}
		In 1962 Charles Hartshorne published a modal logic proof formalizing Anselm of Canterbury's ontological argument for the necessary existence of God. This article presents Kurt G\"odel's notes on this proof which have now been discovered in his \emph{Nachlass} among other theological material, and discusses possible influences on the development of G\"odel's own ontological proof. To complete the picture, strong connections between Anselm of Canterbury's and G\"odel's conceptions of God and his positive properties are pointed out.
	\end{abstract}

\section{Introduction}
	Anselm of Canterbury (1033--1109) is believed to be the first philosopher to have put forward an ontolgical proof or argument for the necessary existence of God. Two variants of this proof can be found both in his little book \emph{Proslogion}, written in 1077/78, as well as in an answer to the monk Gaunilo who criticises the argument as not being valid. In 1962 Charles Hartshorne (1887--2000) published a formal modal logic proof (see \citep{hartshorne}) formalizing Anselm's second version of the argument. While it has been speculated\footnote{See for example \citep{adams}.} whether Kurt G\"odel---when constructing his own ontological proof---was influenced by Hartshorne, he is reported to label the proof as ``wrong'' (see \citep[146]{wang} and section 2), and no direct sign of such an influence could be traced until very recently.\par
	Amongst theological papers in G\"odel's \emph{Nachlass}, we have now discovered (and transcribed from the Gabelsberger shorthand) a single sheet of paper which does not only show that he knew Hartshorne's proof very well. It also presents a remarkable simplification of the proof, accompanied by remarks closely related to Anselm's first variant of the argument as presented in chapter II of the Proslogion, and reappearing in G\"odel's 1970-version \citep{godel1970} in a kind of type-lifted manner.
	
	In order to complete the picture, we show a direct link between Anselm's and G\"odel's conceptions of positiveness and God, thus underlining that G\"odel's ontological proof has been influenced by other sources than Leibniz.

\section{Charles Hartshorne and the ontological proof}

	In 1971/2, Hao Wang conducted extensive interviews with Kurt G\"odel which later formed the basis for his book \citep{wang1996}. In these conversations, G\"odel described Hartshorne as ``an
	example of a contemporary metaphysician.'' \citep[138]{wang1996} This quote is related without
	any deeper context than the fact that ``G\"odel was in favour of metaphysics and opposed to positivism'' which seems to indicate that G\"odel had a generally favourable outlook on Hartshorne in 1971.
	Supporting this, but also bringing up a critisizm of Hartshorne's work, is a later quote from
	1972: When asked to name some recent philosophers whom he found congenial, Gödel replied
	by naming William Henry Sheldon, Josiah Royce and Charles Hartshorne, also describing Hartshorne
	as in the tradition of the scholastics and that he had to ``add what has been done by the followers of Leibniz.'' \citep[141]{wang1996}
	
	However, G\"odel also expressed in 1972 more direct critizism of what he called ``some great
	philosophers of twenty years ago'', including Hartshorne, by attributing the lack of successors to the phenomenon that these philosophers ``generalize things without any inhibition.'' \citep[145]{wang1996} This may also be a specific reference to Hartshorne's ontological proof as using too general logical
	principles as its foundation.
	
	At a later date (Nov 1972) G\"odel again expresses strong criticism of Harts\-horne and his ontological
	proof. G\"odel then states that Hartshorne's lack of knowledge in mathematical logic leads to a
	negative effect resulting in an ontological argument which is wrong.\\
	
	For easier reference we repeat Hartshorne's formal proof along with his original annotations in figure 1, taken from \citep{hartshorne}. The two premises, based on Anselm's argument in chapter III of the Proslogion, can be found in lines 1 and 7, respectively. It should be mentioned that, although $q$ abbreviates the first order formula $(\exists x)P(x)$, a perfect being exists, the steps are purely propositional. Also note that $a \rightarrow b$ abbreviates strict implication or $N\!\sim\!(a\;\& \sim\! b)$ in Hartshorne's notation, a fact that might have been overlooked by G\"odel as we shall see in section 5.\par
	\begin{figure}[t]
		\begin{tabular}{p{3mm}p{3cm}p{5cm}}
			1. & $q \rightarrow Nq$ & ``Anselm's principle''; perfection could not exist contingently\\
			2. & $Nq \;\vee \sim \!Nq$ & Excluded Middle\\
			3. & $\sim Nq \rightarrow N \!\sim\! Nq$ & Form of Becker's Postulate: modal status is always necessary\\
			4. & $Nq \vee N \!\sim\! Nq$ & Inference from (2, 3)\\
			5. & $N \!\sim\! Nq \rightarrow N \!\sim\! q$ & Inference from (1): the necessary falsity of the consequent implies that of the antecedent (Modal form of modus tollens)\\
			6. & $Nq \vee N \!\sim\! q$ & Inference from (4, 5)\\
			7. & $\sim\! N \!\sim q$ & Intuitive postulate (or conclusion from other theistic arguments): perfection is not impossible\\
			8. & $Nq$ & Inference from (6, 7)\\
			9. & $Nq \rightarrow q$ & Modal axiom\\
			10. & $q$ & Inference from (8, 9)
		\end{tabular}
		\caption{Charles Hartshorne's proof.}
	\end{figure}	
	Underneath his proof, Hartshorne adds: ``Those who challange the Argument should decide which of these 10 items or inferential steps to question. Of course one may reject one or more of the assumptions (1, 3, 7); but reject is one thing, refute or show to be a mere sophistry is another.'' Note that he does not include step 2, the law of the excluded middle, in his list of possible rejections.

\section{Simplifications of the proof}
	
	In \citep{schanno2017}, Alisha Schanno presents a simplified version of the above proof, constructed by a theorem prover consisting of no more than 10 lines of \textsc{Prolog} code. As non-logical axioms $Mq$ (\texttt{ip}) and $N(q \rightarrow Nq)$ (\texttt{ap}) are used, which correspond to Hartshorne's lines 7 and 1, respectively. The logical axioms include the modal axioms $MNa \rightarrow Na$ (\texttt{5}), $Na \rightarrow a$ (\texttt{t}) as well as $N(a \rightarrow b) \rightarrow Ma \rightarrow Mb$ (\texttt{km}). The automatically constructed proof of $q$ is then presented as\\[2mm]
	\texttt{mp(t, mp(5, mp(mp(km, ap), ip)))}\\[2mm]
	where \texttt{mp} indicates an application of \emph{modus ponens} or, more correctly, of C.~Mere\-dith's \emph{rule D} of \emph{condensed detachment} which is implemented within the \textsc{Prolog} system as \emph{Robinson's unification algorithm}.\par
	Apart from being much shorter than Hartshorne's version, this proof shows two significant advantages.
	\begin{itemize}
		\item While Hartshorne heavily relies on the law of the excluded middle, Schanno's version is completely intuitionistic, even from a modal perspective.\footnote{The use of the rule \texttt{km} appears to be problematic at first sight, as it can only deduced intuitionistically if the diamond operator $M$ is used as an abbreviation for $\neg N \neg$. (Note that Hartshorne completely avoids the use of the operator $M$.) If both $N$ and $M$ are taken as primitive operators, the above rule has to be added to the underlying logic as an axiom, as was done by Schanno in her proof system. (For further details see \citep{simpson1994}.)}
		\item It follows Eder and Ramharter's criteria \citep{eder2015} for the formal reconstruction of Anselm's ontological argument very closely. Here, the third criterion runs as follows: ``The structure of the formal reconstruction should represent the fundamental structure of the argument. It should be no more and no less detailed than is necessary to map the argument.''
	\end{itemize}
	In order to clarify the second item we reproduce the following table which is taken from \citep{sobel2}, using slightly different notation.\\[4mm]
	
	\begin{center}
		\begin{tabular}{|c|c|c|c|c|}
			\hline 
			& actual world & some poss. world & all poss. worlds & remark\\ 
			\hline 
			1. & $Mq$ &  & & (\texttt{ip}) \\ 
			\hline 
			2. & $N(q \rightarrow Nq)$ &  & & (\texttt{ap})\\ 
			\hline 
			3. &  & $q$ & & (1.)\\ 
			\hline 
			4. & $q \rightarrow Nq$ & $q \rightarrow Nq$ & $q \rightarrow Nq$ & (2.)\\ 
			\hline 
			5. &  & $Nq$ &  & MP (3./4.)\\ 
			\hline 
			6. & $q$ & $q$ & $q$ & (5.)\\ 
			\hline 
		\end{tabular}\\
		Table 1
	\end{center}


	Note how steps 3. to 5. are simulated by axiom \texttt{km}, while the very last step corresponds to an application of the modal axiom 5.
	
	The use of the logic \textbf{S5}\footnote{or related theories involving axiom 5.} in G\"odel's and Hartshorne's proofs has often been criticized as too strong\footnote{For a discussion see \citep{adams}.} and G\"odel himself is reported to have had ``reservations about his ontological proof because of his doubt about using some principle in modal logic.'' \citep[391, fn. \emph{g}]{adams} Whereas modern theorem provers can show that the logic \textbf{KB} is sufficient for the proof (see \citep{benzmuller2014}, \citep{benzmuller2018}), they do not give much insight \emph{why} this is so. The very last line in the above table now clearly reveals the underlying reason: If the necessity of $q$ is established in some possible world, $q$ holds in every possible world, incuding the actual one, $MNq \rightarrow Nq$.\footnote{To make sure that the actual world is included, the axiom T expressing reflexivity is also needed.} But this fact actually carries far too much information for our special case if we are trying to establish God's existence in the \emph{actual} world only. For this purpose it suffices to infer $q$ from the fact that the necessity of $q$ has been confirmed in some possible world, i.e. $MNq \rightarrow q$ (\texttt{b}), thus reducing the underlying logic to \textbf{KB} and the corresponding proof (using Schanno's notation) to\\[2mm]
	\texttt{mp(b, mp(mp(km, ap), ip))}.\\[2mm]
	
	One further reduction of the underlying logic can be based on Charles Hartshorne's own thoughts: He took God's necessary existence as independent of any contingent facts: ``[I]f God logically could be necessary He must be, since no contingent condition can be relevant.'' \citep[53]{hartshorne} Due to God's nature, the proposition $MNq \rightarrow Nq$ could therefore be taken as an additional axiom, thus avoiding the modal axioms T, B and 5 altogether.\footnote{Obviously, Hartshorne did not notice this simplification. This may be due to the fact that his proof uses an equivalent of the modal axiom 5, instantiated by $\sim\! Nq \rightarrow N\!\sim\! Nq$ and annotated as ``Form of Becker's Postulate: modal status is always necessary.''} From a philosophical point of view, this idea may be supported by the question why an underlying (theological) model should necessarily be symmetric or even euclidian as a whole, if only the property $G$ is touched by these structural restrictions. As will be seen in section 6, G\"odel himself would probably have rejected the inclusion of this additional axiom based on the nature of God.

\section{Proslogion II vs. III}
	
	Before we turn to Gödel's own notes on Hartshorne's proof it is worth mentioning Anselm's main argument as presented in Proslogion, chapter II. Here Anselm defines God as ``that than which nothing greater can be conceived.'' His main line of argumentation then runs as follows: As the ``fool'' is able to understand the above definition, God exists in his understanding. Hartshorne formalizes this as $Mq$ (i.e. $\sim\! N \!\sim\! q$ in his terms). But then it is not possible that God exists in the understanding alone because it is greater to exist in the understanding and in reality at the same time. This may be captured as $\sim\! M(Mq \;\& \sim\! q)$ or equivalently as
	\begin{equation}
	 	N(Mq \rightarrow q)
	\end{equation}
	
	In chapter III, Anselm then changes his way of argumentation. Here the central claim is that it is impossible that God exists and at the same time can be thought not to exist, which might be captured as $\sim\! M(q \wedge M\!\sim\! q)$. While Hartshorne repeatedly insisted that the argumentation in chapter III written as
	\begin{equation}
		N(q \rightarrow Nq)
	\end{equation}
	is by far stronger (see for example \citep[11--12, 89]{hartshorne1965}), H.J. Sobel \citep{sobel2} renders the two as equivalent. It can be proved that both modal principles are indeed equivalent in the logic \textbf{KB}.
\section{G\"odel's notes}
	
	G\"odel's notes on Hartshorne's proof which have been discovered in his \emph{Nachlass}\footnote{Kurt G\"odel Papers, Box 10b, Folder 49, item accessions 050156, on deposit with the Manuscripts Division, Department of Rare Books and Special Collections, Princeton University Library. Used with permission of Institute for Advanced Study. Unpublished Copyright Institute for Advanced Study. All rights reserved.} consist of two pages written in logical notation and Gabelsberger shorthand. They are given in figures 2 and 3.\footnote{To improve readability, the overall layout has slightly been altered. Longhand is presented in \emph{italics}.}

	\begin{figure}[t]
		$p \supset Np$ $^\times$ (1)\\
		\uline{$Mp$ $^\circ$ { }}\\
		$Np$ (also $p$)\\[2mm]
		\emph{Allg.} richtiger \emph{Mod. Schluss}\\[4mm]
		\uline{\emph{Bew.}}:\\
		$NNp \vee \sim MNp$,\\
		also $MNp \supset Np$.\\
		Aber aus (1): $Mp \supset MNp$,\\
		also $Mp \supset Np$\\[4mm]
		\uline{\emph{Hartshorne}} \"Uber das \uline{\emph{ontol. Arg.}}\\
		\emph{Logic of Perfect.}, $\overline{62}$, \emph{p.} 50--51\\
		(Was er \uline{in erster Linie zeigt}, ist, dass die Existenz Gottes keine \emph{factual matter} ist.)\\[4mm]
		\uline{$^\circ$ Man muss annehmen, dass man das f\"ur die Existenz Gottes einsehen kann, aber nicht f\"ur die Nicht-Existenz (für seinen Beweis).}\\[2mm]
		$^\times$ Folgt aus dem Wesen Gottes, \uline{also gilt es auch f\"ur die \emph{Negat.}}, aber daraus folgt nicht $\sim p \supset N(\sim p)$. (Z.B. nicht, wenn $p$ wahr und zuf\"allig ist.)\\
		Wohl aber, wenn gilt: $N(p \supset Np)$\\
		Dann gilt auch: $N(\sim p \supset N(\sim p))$\\[2mm]
		\uline{\emph{Bew.}}: Es gilt beides, wenn keine mögliche Welt existiert, in der $p$ (oder keine, in der $\sim p$). Wenn dagegen beide Arten m\"oglicher Welten existieren, dann sind beide falsch.
		\caption{Page 1 of G\"odel's notes. The reference is \citep{hartshorne}.}
	\end{figure}
  
  	\begin{figure}[t]
  			M\"oglichkeit Gottes $\rightarrow$ Existenz Gottes (\emph{Hartsh.}$^\times$)\\[4mm]
  			sehr schlecht kritisiert von \emph{John O. Nelson}\\
  			in \uline{\emph{Rev. of Methaph. 17}, $\overline{63}$, \emph{p. 235--42}}\\[2mm]
  			Antwort darauf: \emph{Rev. of Methaph. 17}, $\overline{63}$, \emph{p. 608}\\[4mm]
  			$^\times$ Er zeigt n\"amlich, dass die Existenz Gottes keine ``factual matter'' ist.
  			\caption{Page 2 of G\"odel's notes. The references are \citep{nelson1963} and \citep{hartshorne1964}, respectively. (Note that Hartshorne's reply was published in 1964.)}
  	\end{figure}
  
  	G\"odel begins his proof with a version of Becker's postulate: modal state is necessary, $NNp \;\vee \sim\! MNp$, either $Np$ holds everywhere or nowhere. He then deduces the modal axiom 5, before he takes the steps presented in \citep{schanno2017}. It is interesting to note that he had already taken these steps in his second version of the ontological proof, presented in \citep{kanckos2019} and written some 10 years earlier.
  	
  	The rest of his thoughts like the footnote marked ($^\times$) written on the first page clearly circulates around the connection between Anselm's two ways of argumentation in Proslogion II and III. Obviously, he takes $p \rightarrow Np$ as a local fact not to be generalized by a rule of necessitation. As mentioned earlier, he might have overlooked the fact that Hartshorne abbreviated (2) by transferring the modal operator $N$ into the implication $\rightarrow$.
	
\section{\emph{Der schlechte Weg} and G\"odel's type-lifting}
	
	In general, G\"odel rejects both Anselm's original form of the argument as well as Hartshorne's modal formalization as can be seen from a note in his notebook \emph{MaxPhil XIV} (p. 107), published in \citep[Appendix B]{godel1995}. Here he states:\footnote{Transcription by the second author, translation taken from \citep[Appendix B]{godel1995}.}
	\begin{quote}
		\small
		\noindent
		Wenn man annimmt $\varphi(x) \supset N\varphi(x)$ [weil aus dem Wesen von $x$ folgend], dann ist es leicht beweisbar, dass es f\"ur jedes kompatible System von Eigenschaften ein Ding gibt, aber das ist der schlechte Weg. Vielmehr soll\footnote{In \citep[Appendix B]{godel1995} wrongly transcribed as ``will''.} $\varphi(x) \supset N\varphi(x)$ erst aus der Existenz Gottes folgen.\par 
		If $\varphi(x) \supset N\varphi(x)$ is assumed [as following from the nature\footnote{In \citep[Appendix B]{godel1995} translated as ``essence''.} of $x$], then it is easily provable that for every compatible system of properties there is a thing, but that is the inferior way. Rather $\varphi(x) \supset N\varphi(x)$ should follow first from the existence of God.
	\end{quote}

	Instantiating $\varphi$ with the predicate $G$ and arguing in the line of Anselm and Hartshorne (``weil aus dem Wesen von $x$ folgend''), \emph{the bad way} (\emph{der schlechte Weg}) now proposes that anything that is God-like is so necessarily. The following ``easy'' (informal) proof starts with this very assumption and immediately leads to G\"odel's conditional claim, where the antecedent expresses the compatibility of the system's (i.e. God's) properties. \par 
	\vspace*{4mm}
	
	\begin{tabular}{p{3mm}p{5cm}p{5cm}}
		(i) & $G(x) \rightarrow NG(x)$ &\\
		(ii) & $\exists x.G(x) \rightarrow \exists x.NG(x)$ &\\
		(iii) & $M\exists x.G(x) \rightarrow M\exists x.NG(x)$ &\\
		(iv) & $M\exists x.G(x) \rightarrow MN\exists x.G(x)$ &\\
		(v) & $M\exists x.G(x) \rightarrow N\exists x.G(x)$ &\\
	\end{tabular} \\[2mm] \par
	
	In \citep{sobel3}, Jordan Howard Sobel takes the above quote as ``strong evidence that G\"odel thought his principles entailed [the modal] collapse, and was not bothered by this, that he indeed considered it a welcome feature of the Leibnizian metaphysics he was devoting for himself.'' However, Sobel misinterprets G\"odel's comment based on a misleading translation. The \emph{Collected Works} \citep[Appendix B]{godel1995} translate ``weil aus dem \emph{Wesen} von $x$ folgend'' as ``as following from the \emph{essence} of $x$'', and Sobel identifies the word ``essence'' with the Leibnizian concept of the thing's ``complete properties'', thus missing the intimate connection to Anselm's argument.\\[2mm]
	
	While clearly rejecting Anselm's argument on the object level,\footnote{See implication (2).} G\"odel adduces a kind of type-lifting of the argument onto the level of properties. The axiom
	\begin{equation}
		P(\varphi) \rightarrow NP(\varphi)
	\end{equation}
	holds, because---following G\"odel's own annotation\footnote{See his 1970 proof \citep{godel1970}. Note the similarity to his ``as following from the nature of $x$.''}---``it follows from the nature of the property.'' Transferring his statement into Anselm's language, \emph{a positive property is so truly positive that it cannot even be thought not to be positive}, thus drawing a parallel to Anselm's argument in \emph{Proslogion} III.
	
	It is worth mentioning that G\"odel's 1970 proof also contains a second variant of the axiom, namely $\neg P(\varphi) \rightarrow N \neg P(\varphi)$. This variant, not being necessary for the proof, did not find its way into Scott's version. Its equivalent
	\begin{equation}
		MP(\varphi) \rightarrow P(\varphi)
	\end{equation}
	again may be taken as a type-lifted version of the argument\footnote{See implication (1).} in \emph{Proslogion} II, \emph{if it is possible that God exists, he exists in reality.} As with Anselm's two ways of argumentation, G\"odel's two variants can quite easily be shown to be equivalent in the logic \textbf{KB}. In the presence of the axiom $P(\varphi) \leftrightarrow \;\sim\! P(\sim \varphi)$, the logic \textbf{K} is already sufficient.\footnote{If the two variants are regarded isolated from the axiom $P(\varphi) \leftrightarrow \;\sim\! P(\sim \varphi)$, a simple countermodel can show that the logic \textbf{K} is not sufficient.} Returning to the aforementioned simplification in section 3, it would be interesting to know if a type-lifted version of Hartshorne's axiom $MNq \rightarrow Nq$, $MN P(\varphi) \rightarrow N P(\varphi)$, could have been included into G\"odel's concept of positiveness, ``because it follows from the nature of the property.''
	
	While it it certainly hard to judge to which extent G\"odel was actually influenced by Hartshorne's proof in this regard, it should be worth mentioning that the axiom $P(\varphi) \rightarrow NP(\varphi)$ appeared in a different shape in exactly those versions of the ontological proof which G\"odel wrote before Harshorne's proof was published in 1962 (see \citep{kanckos2019}). In philosophical notes\footnote{\emph{MaxPhil XIV}, p. 107, published in \citep[434]{godel1995}.} written in 1954 or shortly after he states: ``Da{\ss} die Notwendigkeit einer pos. Eigenschaft pos. ist, ist die wesentliche Voraussetzung f\"ur den ontol. Bew.'' [That the necessity of a positive property is positive is the essential presupposition for the ontological proof.] which is reflected by the axiom $P(\varphi) \rightarrow P(N\varphi)$. And in a footnote relating to the interpretation of positive properties as perfections he mentions the alternative ``Oder: wenn $M\varphi$ eine Perfective, dann auch $\varphi$.'' [Or: if $M\varphi$ is a perfective, then $\varphi$ is, too.] \citep[434]{godel1995}, which in turn may be written as $P(M\varphi) \rightarrow P(\varphi)$.

\section{Essence and existence}
	
	While necessary existence already played a central role even in G\"odel's earliest versions of the ontological proof, the notion of essence only enters the stage in an undated third version, first published in \citep[430]{godel1995} and taken as a complete axiomatization in \citep{kanckos2019}. Here the symbol $ess_x$ playes the role of both a (unique) first order predicate (\emph{the} essence of the object $x$) as well as a second order predicate qualifying a first order predicate as an essence of $x$. While in the first case $ess$ seems to represent a function from objects to predicates, in his 1970-version of the proof $ess$ had become a relation between objects and predicates, thus giving up the (formal) uniqueness of the concept of an essence.
	
	In 1965 Charles Hartshorne expresses a very clear conception of the interplay between the notions of essence and existence, and the resemblance with G\"odel's 1970-version, where he interweaves the two concepts, is striking. In \citep{hartshorne1965} Hartshorne writes:
	\begin{quotation}
		\noindent
		\small
		An important point of the analysis is the inadequacy of the dichotomy \emph{essence-existence}. A third term is needed, 'actuality'. An essence exists if there is some concrete reality exemplifying it; 'existence' is only \emph{that} an essence is concretized, 'actuality' is \emph{how}, or in what particular form, it is concretized. The particular form, the actuality, is always contingent [...] but it does not follow [...] that the existence is contingent. For existence only requires the nonemptiness of the appropriate class of actualities, and a class can be necessarily nonempty even though it has only contingent members. [...]
		
		Essence, existence, actuality---this trias is the minimum of complexity which must be considered if the famous Proof is to be correctly evaluated.
	\end{quotation}
	If existence requires the nonemptiness of the class of actualities, clearly necessary existence requires this very class to have every possible actuality as its member. It is only a small step that leads from necessary existence of an essence to G\"odel's necessary existence of an object via its own essence. It might be speculated whether Hartshorne's existence of an essence (independent of any instantiating object) lead to the much discussed\footnote{See for instance \citep{benzmuller2016} or \citep{benzmuller2018}.} omitting of the conjunct $\varphi(x)$ in G\"odel's definition of essence in 1970, which appears in the preceding version and was later re-added by Dana Scott in \citep{scott1970}. 
	
	\section{Anselm's conception of God}
	
	As has been pointed out in \citep{kanckos2019}, G\"odel got in touch with Anselm and his ontological argument as early as 1925, attending a lecture given by Heinrich Gomperz. When Anselm's conception of God is discussed, the discussion is almost certainly based on the famous phrase ``aliquid quo nihil maius cogitari possit.'' (Proslogion, chapter II) In connection with G\"odel's treatment of the argument it should be noted though, that the Proslogion also involves a very subtle conception based on grades of positiveness. At the end of chapter XI Anselm writes:\footnote{Similar passages appear in chapter V as well as in the \emph{responsio} to the monk Gaunilo.}
	\begin{quotation}
		\noindent
		\small
		Sic ergo vere es sensibilis, omnipotens, misericors et impassibilis, quemadmodum vivens, sapiens, bonus, beatus, aeternus, et quidquid melius est esse quam non esse.\\[2mm]
		So, then, thou art truly sensible, omnipotent, compassionate, and passionless, as thou art living, wise, good, blessed, eternal: and whatever it is better to be than not to be. 
	\end{quotation}
	This is perfectly reflected by G\"odel in his philosophical notes (\emph{MaxPhil XIV}, p. 103, published in \citep[432]{godel1995}) as ``Der ontolgische Beweis mu{\ss} auf den Begriff des Wertes ($p$ besser als $\sim\! p$) gegründet werden.'' [The ontological proof has to be grounded on the concept of value ($p$ better than $\sim\! p$).] which is followed by several axioms concerning the value of properties.
	
	Having captured the notion of positiveness, Anselm goes on to describe the concept of God based on positiveness in chapter XXIV:
	\begin{quotation}
		\noindent
		\small
		Si enim singula bona delectabilia sunt, cogita intente quam delectabile sit illud bonum, quod continet iucunditatem omnium bonorum;\\[2mm]
		For, if individual goods are delectable, conceive in earnestness how delectable is that good which contains the pleasantness of all goods;
	\end{quotation}
	And in the following chapter he continues:
	\begin{quotation}
		\noindent
		\small
		Ama unum bonum, in quo sunt omnia bona, et sufficit. Desidera simplex bonum, quod est omne bonum, et satis est.\\[2mm]
		Love the one good in which are all goods, and it sufficeth. Desire the simple good which is every good, and it is enough.
	\end{quotation}
	These passages reveal a clear correspondence to G\"odel's definition of God, $G(x) \equiv (\varphi)\left[ P(\varphi) \rightarrow \varphi(x) \right] $ which appears in his proofs as early as 1941 (see \citep[429]{godel1995} and \citep{kanckos2019} for a corrected version). They even mirror the impredicativity explicitly surfacing with Dana Scott's axiom $P(G)$, being God-like is itself positive.
	
	The importance of accessing the notion of `God' (\emph{Gottesbegriff}) through his properties is further supported by a set of notes\footnote{Kurt Gödel Papers, Box 6b, Folder 65, item accessions 030088, inserted into the notebook \emph{Max. II Zeiteinteilung}.} mentioning for instance Descartes' ontological proof in a series of loose papers entiteled ``Max. u.\"A.'' (Maximen und \"Ahnliche) and dated ``ca. 1941''. Here G\"odel writes:
	\begin{enumerate}
		\small
		\item[31.] Richtiger Beginn nach \emph{Thomas Aqu. = First Truth =} Wahrheit über Gott, also:
			\begin{enumerate}
				\item[a.)] Die einzelnen S\"atze des Glaubenbekenntnisses in ihrem Sinn klären und Folgerungen ziehen (insbesondere der erste Satz).
				\item[b.)] Von den Eigenschaften Gottes ausgehen, insbesondere \emph{ontolog.} Gottesbeweis (\emph{Descartes}).
				\item[c.)] Insbesondere Abhandlungen lesen, die sich mit dem Gottesbegriff befassen (\emph{August. \& Hilarius, De Trinitate, Summa (De Deo u. Christologie), De divinis nominibus}).
			\end{enumerate}
		\item[31.] Right beginning following Thomas Aqu. = First Truth = truth about God, i.e.:
			\begin{enumerate}
				\item[a.)] To clarify the meaning of the sentences of the Creed and to draw conclusions (in particular the first sentence).
				\item[b.)] To start with God's properties, especially ontological proof of God's existence (Descartes).
				\item[c.)] In particular, to read treatises dealing with the notion of God (August. \& Hilarius, De Trinitate, Summa (De Deo and Christologie), De divinis nominibus).
			\end{enumerate}
	\end{enumerate}
	
	\section{Conclusion}
	In his introduction to G\"odel's ontological proof \citep{adams}, R.M. Adams writes:
	\begin{quotation}
		\noindent
		\small
		Among the historic sponsors of the ontological argument, it is not to Anselm or Descartes but to Leibniz that the parentage of G\"odel's proof belongs, as scholars interested in the proof have long recognized.
	\end{quotation}
	Whereas it has already been shown in \citep{kanckos2019} that G\"odel was well aware of Anselm's argument, the material now found in G\"odel's \emph{Nachlass} clearly underlines the whole breadth of the spectrum of influences on his own versions of the proof, including Anselm of Canterbury, Descartes, and Hartshorne. Nevertheless, it will at the same time be well worth considering the vast material on Leibniz' writings that G\"odel left in his \emph{Nachlass}.
	
\section*{Acknowledgements}
	The research for this article is a part of the \textsc{Godeliana} project led by Jan von Plato, to whom we remain grateful for his support. This project has received funding from the European Research Council (ERC) under the European Union's Horizon 2020 research and innovation programme (grant agreement No 787758) and from the Academy of Finland (Decision No. 318066).

\end{document}